\documentclass[conference]{IEEEtran}
\usepackage[]{times}
\usepackage{epsfig}
\usepackage{subfigure}
\usepackage{amsmath}
\usepackage{bm,amssymb}
\usepackage{graphicx}
\usepackage[top=0.5in, bottom=0.75in, left=0.75in, right=0.75in]{geometry}
\usepackage{multirow}
\usepackage{flushend}
\usepackage{url}
\usepackage{color}

\title{Local Control of Reactive Power by Distributed Photovoltaic Generators}

\author{
\authorblockN{Konstantin Turitsyn}
\authorblockA{CNLS \& Theoretical Divison\\
Los Alamos National Lab\\
Los Alamos,\\ NM 87545, USA\\
Email: turitsyn@lanl.gov}
\and
\authorblockN{Petr \v{S}ulc}
\authorblockA{New Mexico Consortium,\\ Los Alamos, \\ NM
87544, USA\\ Czech Technical University\\ in Prague, \\ Czech Republic\\Email: sulcpetr@gmail.com}
\and
\authorblockN{Scott Backhaus}
\authorblockA{Materials, Physics\\ \& Applications Division\\
Los Alamos National Lab\\
Los Alamos,\\ NM 87545, USA\\
Email: backhaus@lanl.gov}
\and
\authorblockN{Michael Chertkov}
\authorblockA{CNLS \& Theoretical Divison\\
Los Alamos National Lab\\
Los Alamos,\\ NM 87545, USA\\
Also with NMC\\
Email: chertkov@lanl.gov}}

\begin{document}
\maketitle
\begin{abstract}
High penetration levels of distributed photovoltaic (PV) generation on an electrical distribution circuit may severely degrade power quality due to voltage sags and swells caused by rapidly varying PV generation during cloud transients coupled with the slow response of existing utility compensation and regulation equipment.  Although not permitted under current standards for interconnection of distributed generation, fast-reacting, VAR-capable PV inverters may provide the necessary reactive power injection or consumption to maintain voltage regulation under difficult transient conditions.  As side benefit, the control of reactive power injection at each PV inverter provides an opportunity and a new tool for distribution utilities to optimize the performance of distribution circuits, e.g. by minimizing thermal losses.  We suggest a local control scheme that dispatches reactive power from each PV inverter based on local instantaneous measurements of the real and reactive components of the consumed power and the real power generated by the PVs.  Using one adjustable parameter per circuit, we balance the requirements on power quality and desire to minimize thermal losses.  Numerical analysis of two exemplary systems, with comparable total PV generation albeit a different spatial distribution, show how to adjust the optimization parameter depending on the goal. Overall, this local scheme shows excellent performance; it's capable of guaranteeing acceptable power quality and achieving significant saving in thermal losses in various situations even when the renewable generation in excess of the circuit own load, i.e. feeding power back to the higher-level system.
\end{abstract}


\section{Introduction}
\label{sec:intro}
The integration of significant quantities of renewable generation into the U.S. electrical grids presents many significant challenges \cite{lopes2007integrating}.  Reaching this goal will likely require a portfolio of renewable resources and generating methods, e.g. transmission-scale wind, concentrating solar power, and photovoltaic (PV) representing large scale options, while at the small or distribution scale, PV is presently the only viable option. Each of these different generation technologies constitutes a challenge which varies according to the location of its interconnection.  At the transmission scale, renewable generation projects are generally large enough to warrant a transmission interconnection study.  During the study, the transmission operator may uncover problems caused by the renewable generation and require the renewable generation owner to install certain equipment to mitigate the problems.  In this case, the cost of mitigation is borne by the generator creating the transmission problem.

At the distribution scale, the size of an individual PV generator is so small that the cost of an ``interconnection study'' would be prohibitive.  When the number of PV generators on a distribution circuit is small, the impact is quite small and the present utility systems are unaffected.  However, as the PV penetration on a distribution circuit grows, the net impact of many small PV generators can reach a level where the power quality is significantly affected, e.g. loss voltage regulation due rapid variations in PV generation caused by cloud transients coupled with the slow response of existing utility equipment.  Utilities could choose to install fast-response equipment to rectify the problem (e.g. a D-STATCOM \cite{moreno2007power}), however, the cost then would be socialized over the entire rate base and not borne by the PV generators creating the problem and benefiting economically from interconnection to the distribution grid.

An alternate solution is to place the burden of providing reactive power for voltage regulation on the individual, small-scale PV generators by using excess PV inverter capacity to generate or consume reactive power.  However, purposeful injection of reactive power or attempting to regulate voltage by a distributed generator is not permitted by current interconnection standards \cite{1547}.  Changes to these standards to allow for controlling reactive power appear eminent, but even if these changes are made, there still remains significant questions as to how much reactive power to dispatch from each inverter, when to dispatch it, and where and how the control signals are generated.

In contrast to the transmission system, the number and diversity of distribution circuits is perhaps too large to model the individual placement of PV generators.  Therefore, to assess the effectiveness of different algorithms for dispatching reactive power from PV inverters, we take a statistical sampling approach.  We fix certain average `macroscopic' parameters of the distribution circuit: spacing between loads, real and reactive power drawn per load, and PV generation.  However, when creating a realization of the circuit, the individual characteristics of links and values of loads are selected at random from assumed distributions.  For a given set of macroscopic parameters, we create sample realizations of load and generation profiles and analyze the circuit response under the reactive power dispatch algorithm. The macroscopic parameters are varied to assess many different scenarios, e.g. when there are many small PV generators on the circuit but their output is low compared to the circuit loading and when there are only a few PV panels on the feeder line but they inject a substantial amount of power in the circuit leading a net power export and a rise of local voltage.

Analysis such as described above can assess the performance of different algorithms over a wide variety of distribution circuits and can provide valuable guidance to regulatory bodies on questions such as how much excess PV inverter capacity is required to ensure voltage regulation on wide array of distribution circuits.

In our recent work, \cite{KostyaMishaPetrScott} we have focused on the comparison of centralized and decentralized (local) approaches to the control of reactive power. We have shown that, for realistic feeder lines, a simple local control technique can achieve almost 80\% of savings in losses when compared to a centralized control based on solving the full optimization problem. In this work, we extend this approach by discussing the effect of a class of local control techniques on both the losses and power quality in the system.  Although our particular reactive power dispatch algorithm may not work in all cases, the main purpose of this manuscript is to develop a framework for analyzing the effectiveness of any algorithm.

There are other approaches to optimizing the dispatch of reactive power for the purpose of voltage regulation and loss minimization that could be adapted to the present problem including work by Baran and Wu\cite{89BWa,89BWb,89BWc} and Baldick and Wu \cite{90BW} and also in \cite{yona2008optimal,tani2006coordinated}.  However, these works are somewhat specialized to optimal placement and/or control a few large sources of reactive power where the problem at hand includes many small sources.

The remainder of this manuscript is organized as follows.  Section~\ref{sec:Inverter} describes the capability of an inverter to inject or consume reactive power. Section~\ref{sec:DistFlow} describes the power flow solution method we adopt and introduces our reactive power dispatch algorithm.  Section~\ref{sec:rural} provides a description of the prototypical distribution circuit where we test our algorithm.  Section~\ref{sec:simulations} presents the results of our simulations and shows under what conditions our reactive power dispatch algorithm performs well, and Section~\ref{sec:Con} gives our conclusions and suggestions for future work.

\section{Inverter as a limited regulator of local reactive power flow.}
\label{sec:Inverter}

\begin{figure}
\includegraphics[width=0.5\textwidth]{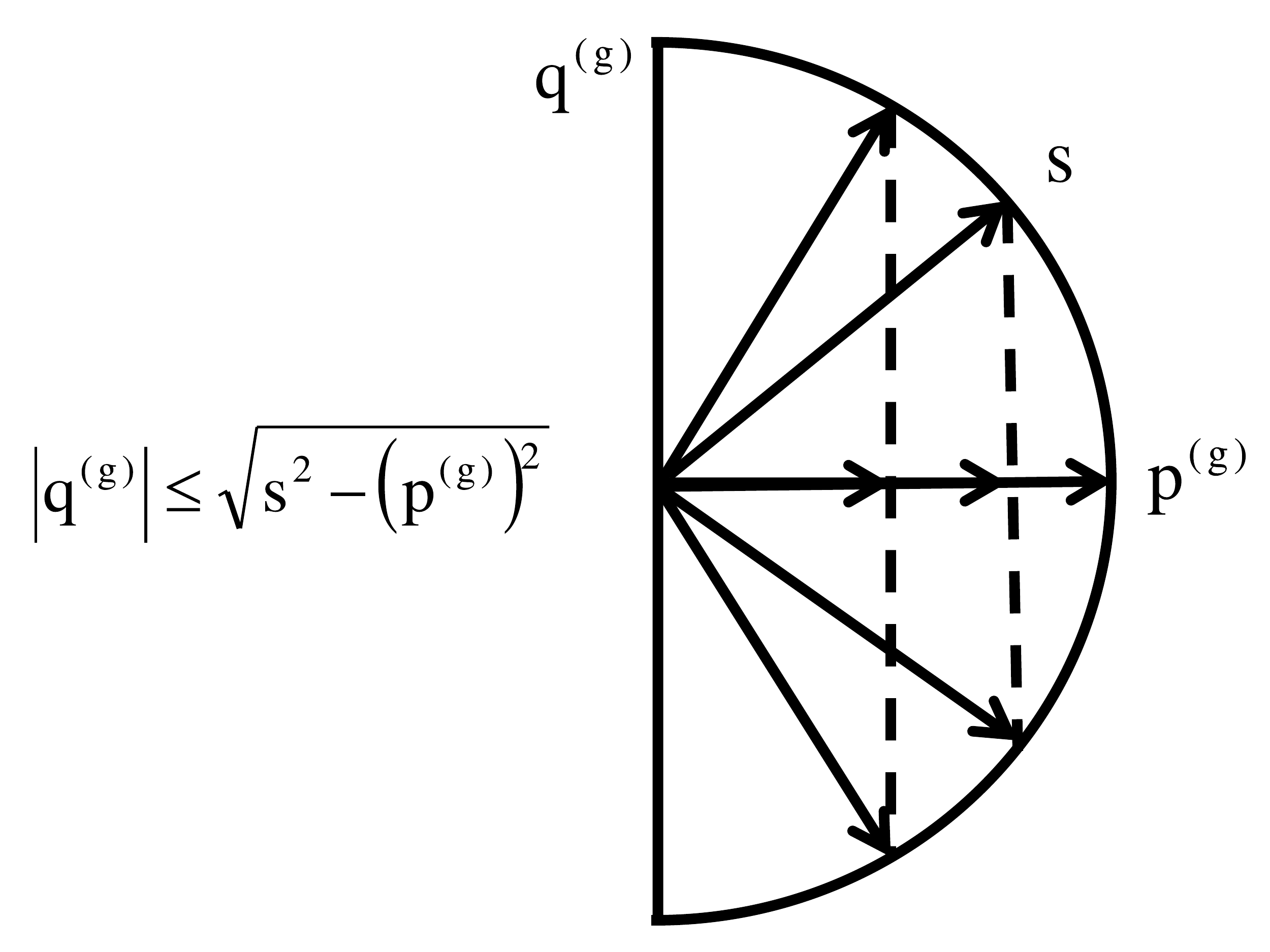}
\centering \caption{When $s$ is larger than $p^{(g)}$, the inverter can supply or consume reactive power $q^{(g)}$.  The inverter can dispatch $q^{(g)}$ quickly (on the cycle-to-cycle time scale) providing a mechanism for rapid voltage regulation.  As the output of the PV panel array $p^{(g)}$ approaches $s$, the range of available $q^{(g)}$ decreases to zero.
}
\label{fig:complex}
\end{figure}

Under the interconnection standard for distributed generation \cite{1547}, PV inverters may not inject or consume reactive power or attempt to regulate voltage in any way, i.e. they must operate at unity power factor while in grid-tied mode.  To overcome voltage regulation problems anticipated on distribution circuits with high penetrations of PV generation, these regulations are expected to be changed to allow the injection of reactive power.  Before we discuss how to dispatch reactive power and analyze the effectiveness of the dispatch, we first discuss the limitations on a PV inverter's reactive power capability.

We adopt a model of PV inverters previously described in \cite{KostyaMishaPetrScott,08LB}.  If the apparent power capability $s$ of an inverter exceed the instantaneous real power generated $p^{(g)}$ by the connected PV panels, the range of allowable reactive power generation is given by $|q^{(g)}|\leq \sqrt{s^2-(p^{(g)})^2}\equiv q^{max}$.  This relationship is also described by the phasor diagram in Fig. \ref{fig:complex}.  On a clear day with the sun angle aligned with the PV array, $p^{(g)}=p^{(g)}_{max}$ and the range of available $q^{(g)}$ is at a minimum.
One focus of this manuscript is to provide a framework for statistically analyzing distribution circuits so set a minimum $s$ relative to $p^{(g)}_{max}$.  However, the analysis is complicated because during cloud events or when the sun angle is not perfectly aligned with the PV panels, $p^{(g)}<p^{(g)}_{max}$ implying that the range of reactive power capability varies though the year, the day, and with the weather.  In order to model this effect we perform our simulations for different values of $p^{(g)}$ while keeping the the absolute capability $s$ fixed.

\section{Optimization of Losses and Voltage Control}
\label{sec:DistFlow}
\begin{figure}
\includegraphics[width=0.5\textwidth]{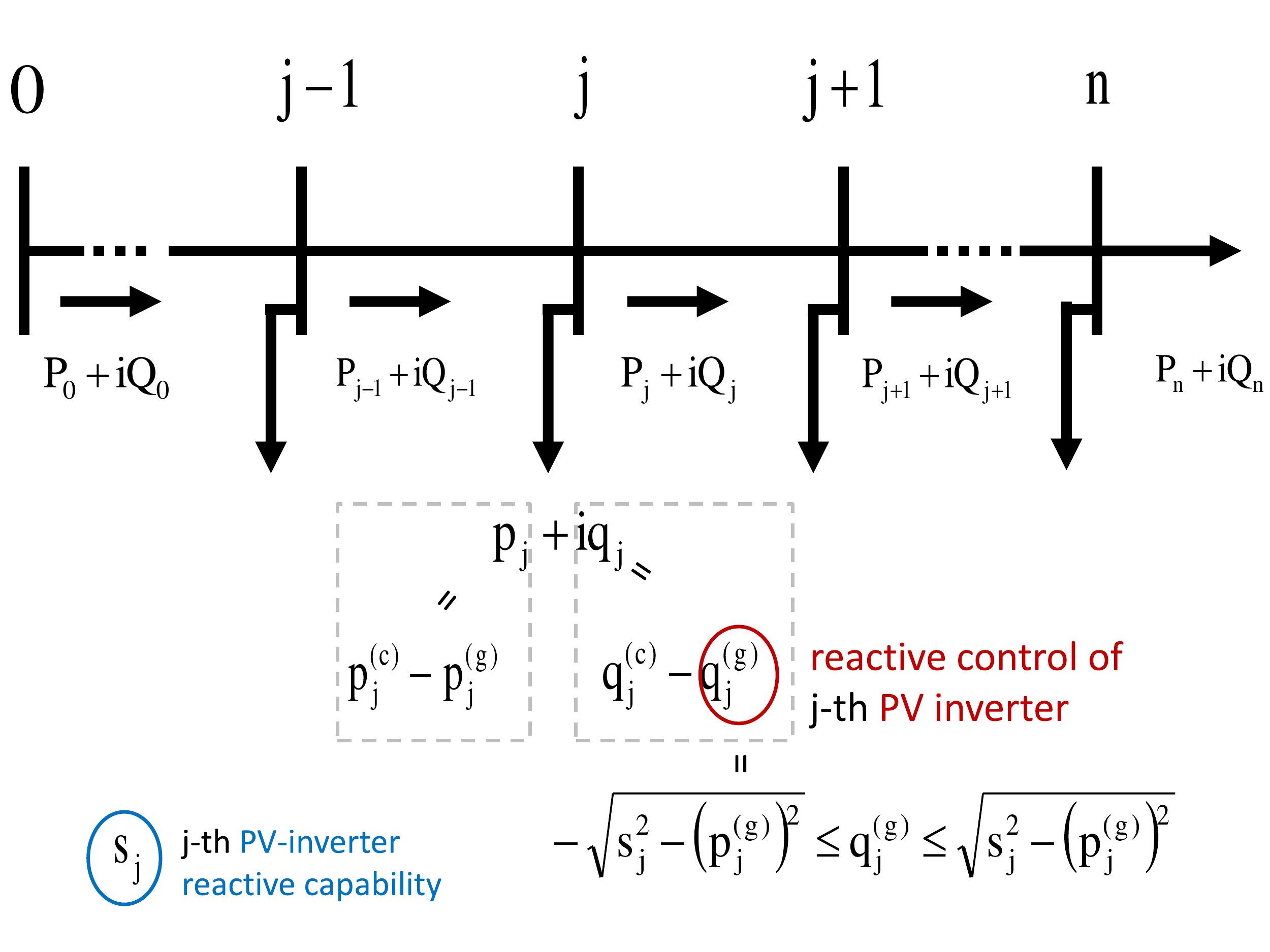}
\centering
\caption{
Diagram and notations for the radial network. $P_j$ and $Q_j$ represent real and reactive power flowing down the circuit from node $j$, where $P_0$ and $Q_0$ represent the power flow from the sub-station. $p_j$ and $q_j$ correspond to the flow of power out of the network at the node $j$, where the respective positive [negative] contributions, $p_j^{(c)}$ and $q_j^{(c)}$ [$p_j^{(g)}$ and $q_j^{(g)}$] represent consumption [generation] of power at the node. The node-local control parameter $q_j^{(g)}$ can be positive or negative but is bounded in absolute value as described in Eq.~\ref{PV_constraint}. The apparent power capability of the inverter $s_j$ is preset to a value comparable to but larger than $\max p_j^{(g)}$ .}
\label{fig:feeder}
\end{figure}

To solve for complex power flows in a radial circuit, we adopt the DistFlow AC power flow equations of \cite{89BWa,89BWb,89BWc}, which for the radial circuit illustrated in Fig.~(\ref{fig:feeder}) can be written $\forall j=1,\cdots,n$:
\begin{eqnarray}
&&P_{j+1}\!=\! P_j\!-\!r_j\frac{P_j^2\!+\!Q_j^2}{V_j^2}\!-\!p_{j+1}, \label{Pj+1}\\
&&Q_{j+1}\!=\!Q_j\!-\!x_j\frac{P_j^2\!+\!Q_j^2}{V_j^2}\!-\!q_{j+1}, \label{Qj+1}\\
&&V_{j+1}^2\!=\!V_j^2\!-\!2(r_jP_j\!+\!x_jQ_j)\!+\!(r_j^2\!+\!x_j^2)
\frac{P_j^2\!+\!Q_j^2}{V_j^2},
\label{Vj2}
\end{eqnarray}
where $P_j+iQ_j$ is the complex power flowing away from node $j$ toward node $j+1$, $V_j$ is the voltage at node $j$, $r_j+ix_j$ is the complex impedance of the link between node $j$ and $j+1$, and $p_j+i q_j$ is the complex power extracted at the node $j$. Both $p_j$ and $q_j$ are composed of local consumption minus local generation due to the PV inverter, i.e. $p_j=p_j^{(c)}-p_j^{(g)}$ and $q_j=q_j^{(c)}-q_j^{(g)}$.  Of the four contributions to $p_j+i q_j$, $p_j^{(c)}$, $p_j^{(g)}$ and, $q_j^{(c)}$ are uncontrolled (i.e. driven by consumer load or instantaneous PV generation), while the reactive power generated by the PV inverter, $q_j^{(g)}$, can be adjusted.  As described in Section~\ref{sec:Inverter}, $q_j^{(g)}$ is limited by the reactive capability of the inverter:
\begin{eqnarray}
\forall j=1,\cdots,n:\quad
\left| q_j^{(g)} \right| \leq \sqrt{s_j^2-(p_j^{(g)})^2} \equiv q_j^{max}. \label{PV_constraint}
\end{eqnarray}
Note that reactive power generation is possible only at the nodes with PV generation. For the other nodes, we take $s_j=0$.
For a broad class of feeder lines such as the one considered in this work, the quadratic terms in Eqs.~(\ref{Pj+1},\ref{Qj+1},\ref{Vj2}) are relatively small as are the deviations of voltage along the line. In this case the power flow equations can be approximated by the linear \emph{LinDistFlow} equations \cite{89BWa,89BWb,89BWc}:
\begin{eqnarray}
\label{LinDistFlow1}
P_{j+1}= P_j-p_{j+1}^{(c)}+p_{j+1}^{(g)},\\
Q_{j+1}=Q_j-q_{j+1}^{(c)}+q_{j+1}^{(g)}, \label{LinDistFlow2}\\
V_{j+1}=V_j - (r_j P_j+x_j Q_j)/V_0. \label{LinDistFlow3}
\end{eqnarray}
where we have exploited the approximation $V_k^2 \approx V_0^2 + 2 V_0 \left( V_k - V_0 \right) $.

Within the framework of linearized model, the rate of energy dissipation (losses) in the distribution circuit is given by
\begin{eqnarray}
{\cal L}=\sum_{j=0}^{n-1} r_j\frac{P_j^2+Q_j^2}{V_0^2},
\label{Loss}
\end{eqnarray}
Minimizing or at least keeping the losses acceptably low is a natural goal for optimization and control. However, voltage variations along the circuit must stay within strict regulation bounds.  These bounds are respected if the maximum deviation of the per unit voltage $\delta V$ obeys
\begin{equation}\label{Vdef}
  \delta V = \max_{k} \left| \frac{V_k-V_0}{V_0} \right|<\epsilon,
\end{equation}
where $\epsilon\approx 0.05$ in normal operation. The goal of the voltage stability optimization is to keep $\delta V$ within normal bounds.

In this work, we discuss local control techniques that try to achieve objectives of minimizing both the thermal losses (\ref{Loss}) and the maximal voltage deviation (\ref{Vdef}).  Our local scheme should be viewed as an approximate/heuristic solution of the following multi-objective optimization problem
\begin{eqnarray}
\label{LinDistFlow}
&&\min_{{\bm q}^{(g)}} \left[{\cal L}, \delta V \right]^T,\\
&&\mbox{s.t.} \,\mbox{Eqs.~(\ref{PV_constraint},\ref{LinDistFlow1},\ref{LinDistFlow2},\ref{LinDistFlow3})}.\nonumber
\end{eqnarray}
The locality of the control scheme corresponds to the requirement that the $q_k^{(g)}$ depend only on local information, which we restrict to $p_k^{(g)}, p_k^{(c)}, q_k^{(c)}$. Another local variable that could be used is $V_k$, however, using $V_k$ could easily lead to inequities.  For instance, compared to a PV generator near the beginning of a distribution circuit, one at the end of the circuit may perpetually see a relatively lower voltage and be requested to inject reactive power leading to faster degradation of the inverter.  Instead, our scheme is based on local real and reactive power consumption and generation and should be more equitable along the circuit.  We also assume that the control scheme is homogeneous over the line: all inverters are programmed in the same way, and explicit dependence on the bus number $k$ enters through the inverter capability $s_k$. Formally, we will study local control schemes with
\begin{equation}\label{control}
	q_k^{(g)} = F_k(p_k^{(g)}, p_k^{(c)}, q_k^{(c)}),
\end{equation}
and consistent with constraint (\ref{PV_constraint}). It is useful to introduce the following
``helper" function, $\text{Constr}_{j}$, meant to enforce the constraint \eqref{PV_constraint}:
\begin{equation}\label{multi}
 \text{Constr}_{j}[q] = \begin{cases}
q, &    \left| q \right|  \leq  q_j^{max}  \\
(q/|q|)q_j^{max}, & \text{otherwise}
\end{cases}
\end{equation}
The local control scheme from \cite{KostyaMishaPetrScott} was designed based on the idea that losses are minimized when the reactive flows $Q_k$ are zero. The $q_k^{(g)}$ were chosen to minimize the {\it net} reactive power consumption $q_k^{(c)}-q_k^{(g)}$ at each node:
\begin{equation} \label{losscontrol}
	F_k^{(L)} = \text{Constr}_{k}[q_k^{(c)}].
\end{equation}
This scheme was shown to be very effective in reducing the losses. In this work, we extend this idea and propose a new class of the schemes that attempt to achieve simultaneously both optimization objectives in Eq.~(\ref{LinDistFlow}).

Eq.~(\ref{LinDistFlow3}) suggests that, to reduce variations in $V_k$, we should minimize the absolute value of the combined power flow $r_k P_k + x_k Q_k$. Note that for many circuits, the ratio of $r_k/x_k = \alpha$ is nearly constant for all $k$ and depends only on the configuration and size of the conductors used. Thus the absolute value of $r_k P_k + x_k Q_k$ will be exactly zero if for every line we will ensure that $p_j^{(c)} - p_j^{(g)} + \alpha \left( q_j^{(c)} - q_j^{(g)} \right) = 0$. Summarizing,  this suggests the following control function $F_k^{(V)}$, aimed at minimizing  voltage variations and ignoring losses:
\begin{equation}\label{voltcontrol}
	F_k^{(V)} = \text{Constr}_{k}\left[  q_k^{(c)} + \frac{p_k^{(c)} - p_k^{(g)}}{\alpha} \right].
\end{equation}
{In essence, $F_k^{(V)}$ is attempting to not only supply the reactive power needs of the local load at node $k$, but also to supply the reactive power consumption of the adjacent links of the distribution circuit.  A compromise between the two objectives in (\ref{multi}) can be achieved via the following nonlinear combination
\begin{equation}\label{hybrid}
	F_k = \text{Constr}_{k}\left[ K F_k^{(L)} + (1-K) F_k^{(V)} \right],
\end{equation}
where $K$ is a single parameter controlling the trade off between the two objectives in Eq.~(\ref{multi}). At $K=1$ we recover the (\ref{losscontrol}) scheme, whereas at $K=0$ the scheme reduces to (\ref{voltcontrol}).

\section{Description of the prototypical rural distribution circuit}
\label{sec:rural}
We consider a sparsely-loaded rural distribution circuit model with $250$
nodes based on one of the 24 prototypical distribution circuits described
in \cite{08SCCPET}. The nominal phase-to-neutral voltage $V_0$ is $7.2 kV$ and the line impedance is set to
$(0.33+0.38i) \Omega /km$, constant along the distribution line. The
distance between two neighboring nodes is uniformly distributed between $0.2$
and $0.3$ kilometers. The real powers $p_j^{(c)}$ consumed by the loads are uniformly distributed between $0$ and $p^{(c,max)}$ where $p^{(c,max)}$ is
set to either to $1 kW$ or $2.5 kW$ (for the different cases we consider), and $q_j^{(c)}$ is randomly selected from uniform distribution between $0.2 p_j^{(c)}$ and $0.3 p_j^{(c)}$. In this
work, we study two levels of PV penetration; either
$20\%$ or $50\%$ of the nodes have PV generation installed, but the PV-enabled
nodes are selected randomly. We assume
that the generated power $p_j^{(g)}$ is the same for all PV-enabled nodes being equal to either $1 kW$ or $2 kW$. A uniform level of PV power generation
assumes identical installations at each PV-enabled node (the same $p^{(g)}_{max}$ installed in the same way) and uniform solar irradiance.  The inverter capacity is set to $s = 2.2 kW$ at PV-enabled nodes.

The algorithm we use for the simulation consists of several steps:
first we generate the random values of $r_k$ and $x_k$ that remain constant for all
forthcoming simulations. For every case, we generate a random
sample of loads and generation and solve the LinDistFlow equations
(\ref{LinDistFlow1}-\ref{LinDistFlow3}) to find both the voltage levels along the line
and the total losses. We do not perform any averaging over many realizations of the load distribution because the results do not differ significantly from one sample to
another due to the self-averaging effect of a large number of node on the line.

\section{Simulations: Results and Discussions}
\label{sec:simulations}

\begin{table}
\caption{Comparison of parameters for different model feeder lines}
\begin{center}
\begin{tabular}{ | c | c | c | c | c | c |}
    \hline
    Case \# & PR & $p^{(c,mean)}$ & $p^{(g)}$ & $\delta V_0$ & ${\cal L}_0$\\ \hline
    1 & 20 \%  & 1.25 kW & 1.0 kW & 0.059 & 7.84 kW \\
    2 & 20 \%  & 0.5 kW & 2.0 kW & 0.014 & 0.33 kW \\
    3 & 50 \%  & 1.25 kW & 1.0 kW & 0.048 & 4.66 kW \\
    4 & 50 \%  & 0.5 kW & 2.0 kW & 0.014 & 1.89 kW \\ \hline
\end{tabular}
\end{center}
\label{tab:cases}
\end{table}

Four regimes of operation are considered corresponding to two different levels of PV penetration and two operational extremes (a) sunny day, everybody is at work, and (b) overcast,  everybody home.  These are summarized in Table \ref{tab:cases}. In cases $1$ and $3$ the amount of power generated is low compared to the total load of the system, however the capacity of the inverters is high: $s/p^{(g)} \approx 2.2$. In the other two cases the situation is opposite; the total amount of power generated is almost equal to the total consumption in case $2$, and there is over generation of power in case $4$ leading to reversal of the power flow (the circuit is feeding power into the higher level grid). The inverter capacity in cases $2$ and $4$ is relatively small: $s/p^{(g)} \approx 1.1$. These $4$ cases allow us to assess the importance of several characteristics of the system: generation and load levels, PV penetration rate, and inverter capacity.

The reactive power is controlled according to Eq.~(\ref{hybrid}) with the parameter $K$ set to the same level on all the generating nodes. $K$ is varied in the range $-5 < K < 10$, and the  total losses ${\cal L}$ and maximal voltage deviation $\delta V$ are computed. For all cases, the performance of the control technique is compared to the base case of zero reactive power generation: $q^{(g)}_k = 0$. The corresponding values of ${\cal L}_0$ and $\delta V_0$ are shown in Table \ref{tab:cases}.

\begin{figure}
\includegraphics[width=0.5\textwidth]{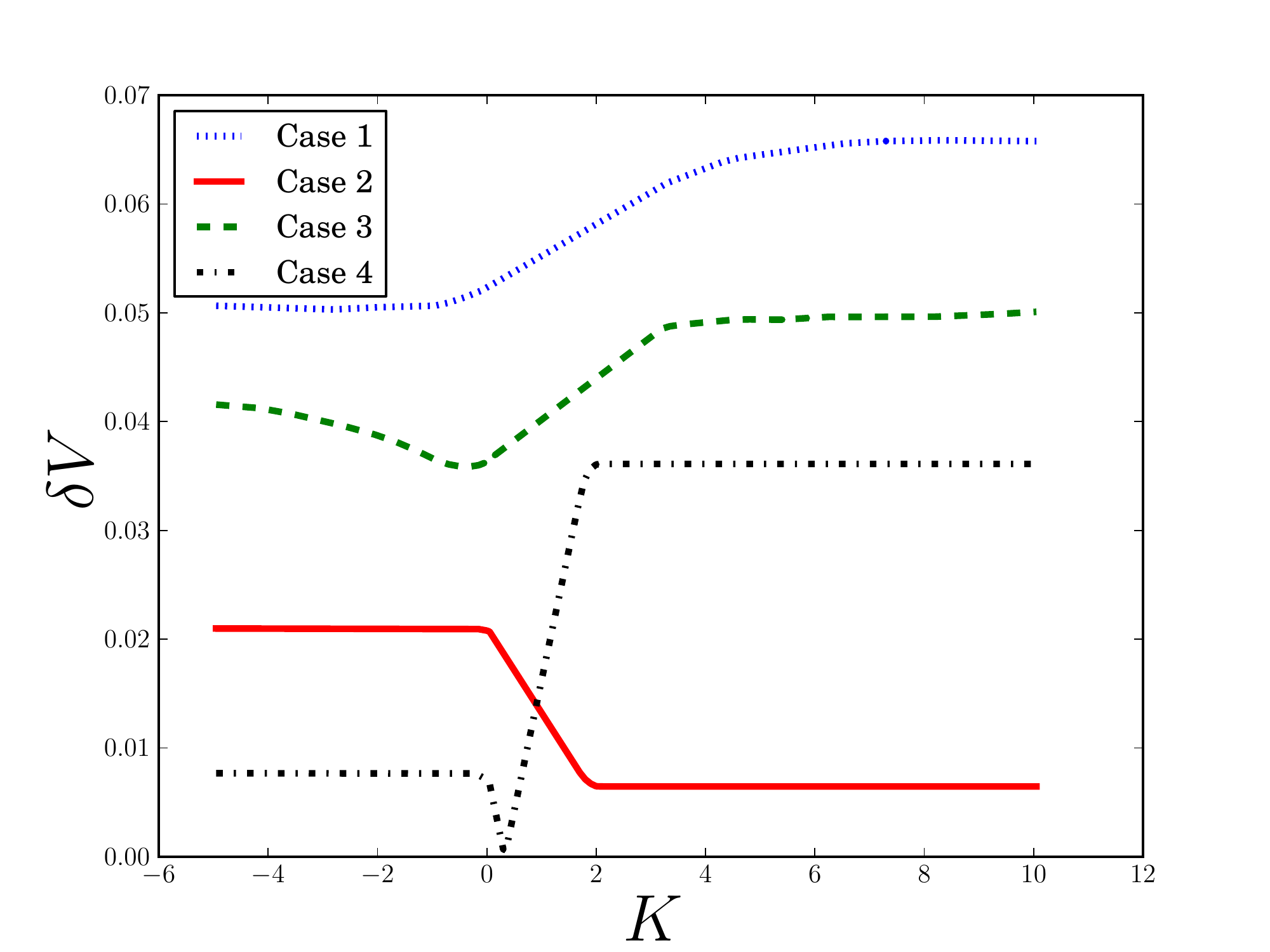}
\centering \caption{
$\delta V  = \max_{k} \left| \frac{V_k-V_0}{V_0} \right|$ vs $K$ for the local control according to \eqref{hybrid} and four cases described in Table \ref{tab:cases}.}
\label{fig:V}
\end{figure}

The dependence $\delta V$ versus $K$ is shown in Fig. \ref{fig:V}. The performance of the voltage control is very good in cases $2$ and $4$, when there is a small amount of power consumption. Significant reduction in voltage deviation in the case $4$ is possible because of almost complete compensation of the reversed flow of real power by increase in the consumption of reactive power. In cases $1$ and $3$, correspondent to high power consumption, the reduction in voltage deviation is less significant, however it is still possible to achieve the reduction of $\approx 0.01$ in both cases.

\begin{figure}
\includegraphics[width=0.5\textwidth]{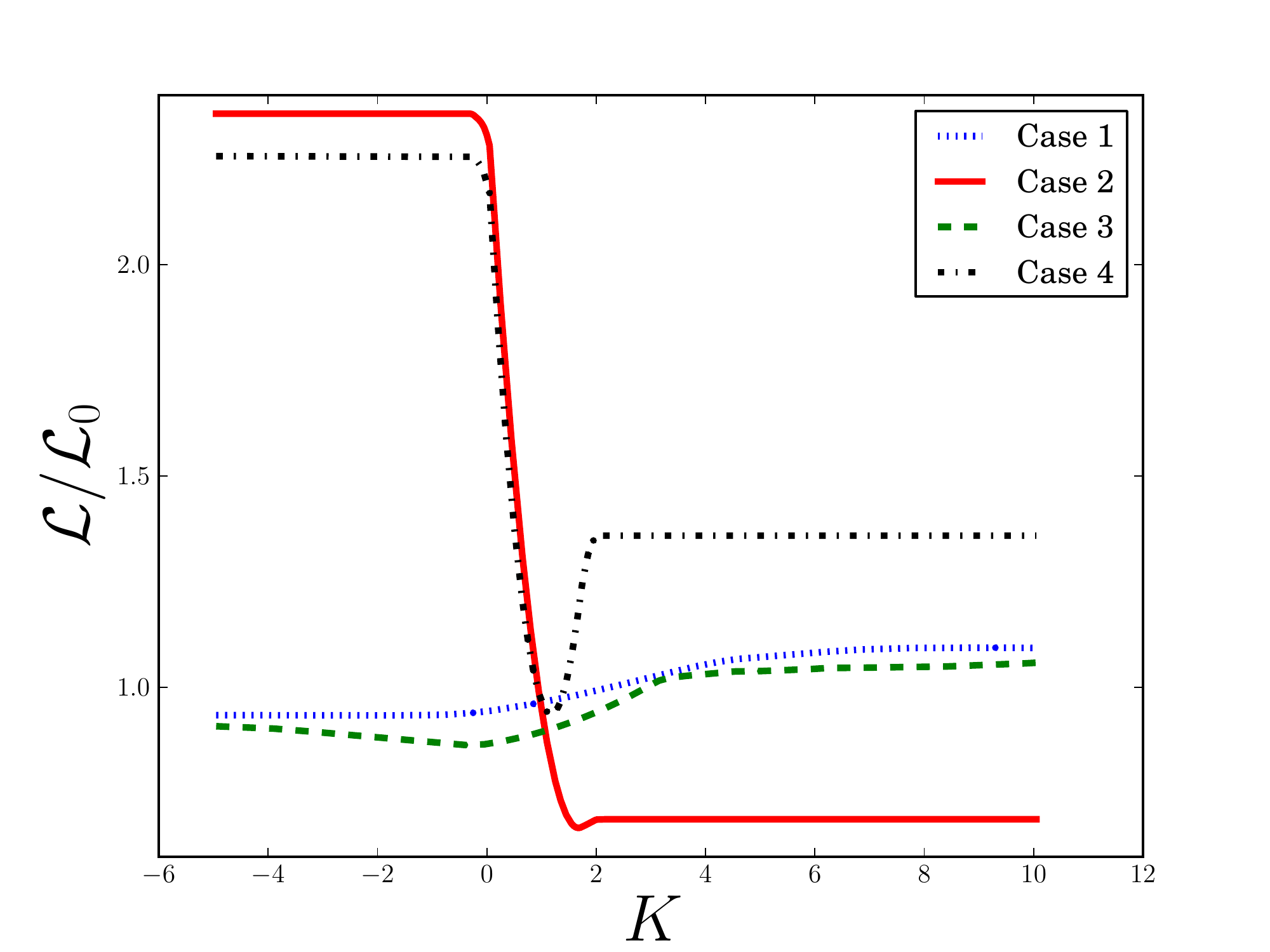}
\centering \caption{Dissipated power with control, measured in the units of the respective dissipated power without control, vs $K$ for the local control according to \eqref{hybrid} and four regimes described in Table \ref{tab:cases}.
The lowest relative losses achievable are $0.93$, $0.67$, $0.86$ and $0.94$ for cases $1$ to $4$ respectively.
 }
\label{fig:losses}
\end{figure}

Changes in losses with variation in the parameter, $K$, shown in Fig.~\ref{fig:losses}, are similar to these of the power quality discussed above. The top savings are achieved in the cases $2$ and $3$. Thus, at the optimal point of the case $2$ the consumption is almost completely balanced by the generation, and then the magnitudes of total power flows $P_k,Q_k$ are greatly reduced. Note, that although the relative savings in losses are very impressive in cases $2$ and $4$, the absolute values are still smaller compared to these in the cases $1$ and $3$ because of the lower values of ${\cal L}_0$ in the latter. The high absolute savings in the cases $1$ and $3$ should be attributed to the high value of the relative inverter capacity: $s/p^{(g)} = 2.2$, that allows high injections of reactive power in the system. Despite, the high penetration level the savings in losses in the case $4$ are less pronounced in comparison to the case $2$ because of significant losses associated with the reverse flow of real power.  Interestingly, and in a contrast to other cases, optimization of losses in the case $4$ required fine tuning of the parameter $K$.

\begin{figure}
\includegraphics[width=0.5\textwidth]{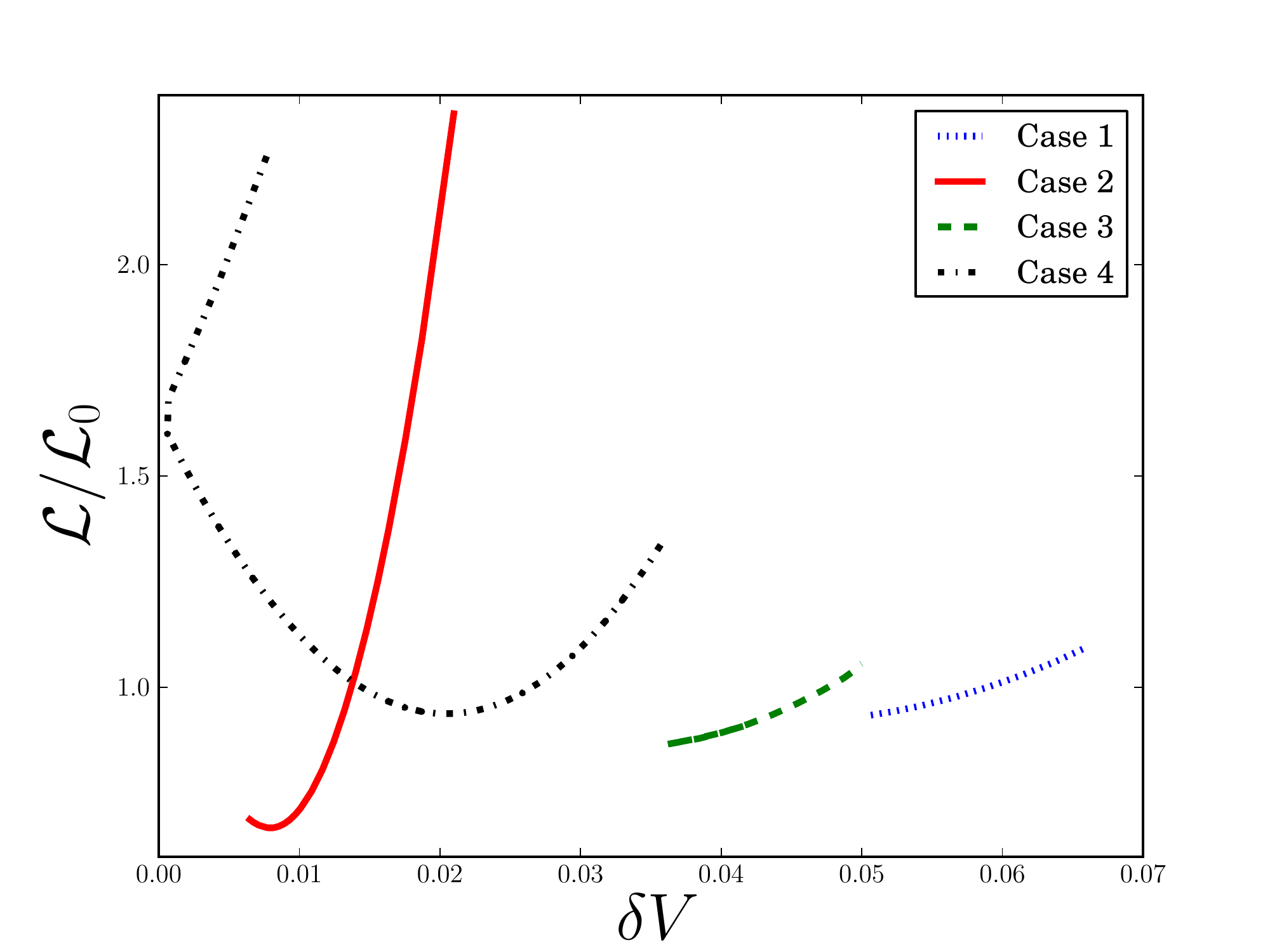}
\centering \caption{The plot combines results from Fig.~\ref{fig:losses} and Fig.~\ref{fig:V} in one, thus showing dependence of the relative power loss vs power quality characteristic,  when scanned over allowed range of the parameter $K$ in control scheme \ref{hybrid}. In the cases 1 and 3, the left most point of the curve corresponds to the best configuration, achieving both the lowest $\delta V$ and the smallest losses. For cases 3 and 4, the whole part of the curve from the left most  point (corresponding to the best power quality) to the point where global minimum is achieved (losses are minimal then) represents a feasible region where one can trade-off between savings in losses and the voltage drop by adjusting $K$.
}
\label{fig:Vloss}
\end{figure}

The overall performance of the control scheme \ref{hybrid} can be most easily assessed on the combined plot on Figure \ref{fig:Vloss} where the efficiency in both $\delta V$ and ${\cal L}$ is shown simultaneously. It is apparent from this plot that in the cases $1$ and $3$ the optimum in both losses and voltage deviations is achieved simultaneously. From the viewpoint of the multi-objective optimization this corresponds to Pareto set consisting of a single point. The situation is different in the cases $2$ and $4$ where the Pareto set is represented on a range in $K$. (This set corresponds to the negative slope parts of the curves in Fig. \ref{fig:Vloss}.)  Comparing different cases, the maximal flexibility is achieved in the case $4$ where the high penetration of renewables results in the overgeneration of power.


\section{Conclusions and Path Forward}
\label{sec:Con}
To summarize, this study suggested
\begin{itemize}
    \item Simple and local control scheme (\ref{hybrid}) for the reactive power generation at circuit-distributed PV inverters.  The control adjusts the reactive power according to the local values of (real and reactive) consumption and real PV generation. The scheme contains one global adjustable parameter balancing between complementary (local) desires to minimize the power flow and to maintain good power quality.
     \item We experimented with the scheme showing that the local scheme is capable of achieving very significant (and probably sufficiently close to globally optimal) simultaneous improvement in global power quality and global reduction of losses over the entire circuit. Different regimes tested correspond to sunny/overcast and large/small load cases.
    \item We concluded that already $10\%$ of excess inverter capacity $s$ is enough to allow significant reductions in both the losses and voltage deviations.
    \item We observe that in the case of the over-generation, and thus reversal of the power flow direction,  the multi-objective optimization becomes more sensitive to the parameter change,  thus requiring a more accurate tuning.
\end{itemize}

The results reported in this study will require further analysis and exploration along the following lines:
\begin{itemize}
    \item Study of other local control schemes aimed at even better improvements.
    \item We conjectured that fluctuations from an instance to instance (e.g. in configurations of load) drawn from a reasonable distribution representing a typical regime (say a morning of the overcast day) will not be significant.  This conjecture needs to be tested in a more accurate statistical study.  The study should test statistical variations in the control scheme performance.
    \item It may be beneficial for improving the power quality characteristics to include information about local voltage (if such information is available) into the control scheme.
    \item We also envision generalizing the scheme to a global dynamical control,  with the parameter $K$ varying in time according to a broadcasted central command. Moreover,
        a distributed control,  with the coefficient $K$ also varying spatially,  may also be considered as beneficial in special (geographically varying) situations,  where we should also consider inhomogeneity and significant variation in the instantaneous insulation pattern.
\end{itemize}

\section*{Acknowledgment}

We are thankful to all the participants of the ``Optimization and Control for Smart Grids" LDRD DR project at Los Alamos
and Smart Grid Seminar Series at CNLS/LANL for multiple fruitful discussions. Research at LANL was carried out under the auspices of the National Nuclear Security Administration of the U.S. Department of Energy at Los Alamos National Laboratory under Contract No. DE C52-06NA25396. P\v{S} and MC acknowledges partial support of NMC via NSF collaborative grant CCF-0829945 on ``Harnessing Statistical Physics for Computing and Communications''.

\bibliographystyle{IEEEtran}
\bibliography{SmartGrid}

\begin{thebibliography}{10}
\providecommand{\url}[1]{#1}
\csname url@samestyle\endcsname
\providecommand{\newblock}{\relax}
\providecommand{\bibinfo}[2]{#2}
\providecommand{\BIBentrySTDinterwordspacing}{\spaceskip=0pt\relax}
\providecommand{\BIBentryALTinterwordstretchfactor}{4}
\providecommand{\BIBentryALTinterwordspacing}{\spaceskip=\fontdimen2\font plus
\BIBentryALTinterwordstretchfactor\fontdimen3\font minus
  \fontdimen4\font\relax}
\providecommand{\BIBforeignlanguage}[2]{{%
\expandafter\ifx\csname l@#1\endcsname\relax
\typeout{** WARNING: IEEEtran.bst: No hyphenation pattern has been}%
\typeout{** loaded for the language `#1'. Using the pattern for}%
\typeout{** the default language instead.}%
\else
\language=\csname l@#1\endcsname
\fi
#2}}
\providecommand{\BIBdecl}{\relax}
\BIBdecl

\bibitem{lopes2007integrating}
J.~Lopes, N.~Hatziargyriou, J.~Mutale, P.~Djapic, and N.~Jenkins,
  ``{Integrating distributed generation into electric power systems: A review
  of drivers, challenges and opportunities},'' \emph{Electric Power Systems
  Research}, vol.~77, no.~9, pp. 1189--1203, 2007.

\bibitem{moreno2007power}
A.~Moreno-Munoz, \emph{{Power quality: mitigation technologies in a distributed
  environment}}.\hskip 1em plus 0.5em minus 0.4em\relax Springer Verlag, 2007.

\bibitem{1547}
\BIBentryALTinterwordspacing
``{IEEE 1547 Standard for Interconnecting Distributed Resources with Electric
  Power Systems}.'' [Online]. Available:
  \url{http://grouper.ieee.org/groups/scc21/1547/1547_index.html}
\BIBentrySTDinterwordspacing

\bibitem{KostyaMishaPetrScott}
K.~Turitsyn, P.~\v{S}ulc, M.~Chertkov, and S.~Backhaus, ``Use of reactive power
  flow for voltage stability control in radial circuit with photovoltaic
  generation,'' in \emph{Power Engineering Society General Meeting, 2010.
  IEEE}, July 2010.

\bibitem{89BWa}
M.~Baran and F.~Wu, ``Optimal sizing of capacitors placed on a radial
  distribution system,'' \emph{Power Delivery, IEEE Transactions on}, vol.~4,
  no.~1, pp. 735--743, Jan 1989.

\bibitem{89BWb}
------, ``Optimal capacitor placement on radial distribution systems,''
  \emph{Power Delivery, IEEE Transactions on}, vol.~4, no.~1, pp. 725--734, Jan
  1989.

\bibitem{89BWc}
------, ``Network reconfiguration in distribution systems for loss reduction
  and load balancing,'' \emph{Power Delivery, IEEE Transactions on}, vol.~4,
  no.~2, pp. 1401--1407, Apr 1989.

\bibitem{90BW}
R.~Baldick and F.~Wu, ``Efficient integer optimization algorithms for optimal
  coordination of capacitors and regulators,'' \emph{Power Systems, IEEE
  Transactions on}, vol.~5, no.~3, pp. 805--812, Aug 1990.

\bibitem{yona2008optimal}
T.~Yona and N.~Funabashi, ``{Optimal Distribution Voltage Control and
  Coordination With Distributed Generation},'' \emph{IEEE Transactions on Power
  Delivery}, vol.~23, no.~2, 2008.

\bibitem{tani2006coordinated}
J.~Tani and R.~Yokoyama, ``{Coordinated Allocation and Control of Voltage
  Regulators Based on Reactive Tabu Search for Distribution System},''
  \emph{WSEAS Transactions on Power Systems}, vol.~1, no.~2, 2006.

\bibitem{08LB}
\BIBentryALTinterwordspacing
E.~Liu and J.~Bebic, ``Distribution system voltage performance analysis for
  high-penetration photovoltaics,'' NREL/SR-581-42298, Tech. Rep., 2008.
  [Online]. Available: \url{http://www1.eere.energy.gov/solar/pdfs/42298.pdf}
\BIBentrySTDinterwordspacing

\bibitem{08SCCPET}
\BIBentryALTinterwordspacing
K.~Schneider, Y.~Chen, D.~Chassin, R.~Pratt, D.~Engel, and S.~Thompson,
  ``Modern grid initiative-distribution taxonomy final report,'' Tech. Rep.,
  2008. [Online]. Available:
  \url{http://www.gridlabd.org/models/feeders/taxonomy_of_prototypical_feeders%
.pdf}
\BIBentrySTDinterwordspacing

\end{thebibliography}

\end{document}